\documentclass[a4paper,10pt,reqno]{amsart}
\usepackage{amssymb,latexsym,bbm,amsmath,amsfonts}
\usepackage{amsthm}
\usepackage[mathscr]{eucal}
\usepackage{rotating}
\usepackage{multirow}
\usepackage{slashbox}

\numberwithin{equation}{section}
\newtheorem{theorem}{Theorem}

\newtheorem{proposition}{Proposition}

\theoremstyle{definition}
\newtheorem{definition}[theorem]{Definition}

\newtheorem{example}[theorem]{Example}

\theoremstyle{remark}
\newtheorem{remark}{Remark}

\newcommand{\norm}[1]{\left\lVert #1 \right\rVert}

\DeclareMathOperator{\R}{{\mathbb R}}                
\DeclareMathOperator{\Rp}{{\mathbb R}_+}             
\newcommand{\e}{\mathrm{e}}                                                  

\begin{document}

\title[Admissibility of Linear Stochastic Volterra Operators]
{On the Admissibility of Linear Stochastic Volterra Operators}

\author{John~A.~D.~Appleby}
\address{Edgeworth Centre for Financial Mathematics, School of Mathematical
Sciences, Dublin City University, Glasnevin, Dublin 9, Ireland}
\email{john.appleby@dcu.ie} \urladdr{webpages.dcu.ie/\textasciitilde
applebyj}

\author{John~A.~Daniels}
\address{Edgeworth Centre for Financial Mathematics, School of Mathematical
Sciences, Dublin City University, Glasnevin, Dublin 9, Ireland}
\email{john.daniels2@mail.dcu.ie}

\author{David~W.~Reynolds}
\address{School of Mathematical
Sciences, Dublin City University, Glasnevin, Dublin 9, Ireland}
\email{david.reynolds@dcu.ie}

\thanks{Both authors gratefully acknowledge Science Foundation Ireland for the support of this research
under the Mathematics Initiative 2007 grant 07/MI/008 ``Edgeworth
Centre for Financial Mathematics''. The second author is
also supported by the Irish Research Council for Science,
Engineering and Technology under the Embark Initiative grant.}
\date{11 October 2012}

\begin{abstract}
Conditions guaranteeing convergence of linear stochastic Volterra
operators are studied. Necessary and sufficient conditions for mean
square convergence are established, while almost sure convergence of
the linear operator is shown to imply mean square convergence.
Sufficient conditions for almost sure convergence of the stochastic
linear operator are established.
The sharpness or necessity of these conditions is explored by means of examples.
\end{abstract}

\maketitle


\section{Introduction}
This article studies convergence properties of It\^o-Volterra integrals of the form
%
\begin{equation}\label{eq:HfB}
(\mathcal{H}f)(t):=\int_0^t H(t,s)f(s)\,dB(s)
\end{equation}
where $H$ is a deterministic Volterra kernel and $f$ is a
deterministic function on $[0,\infty)$. We require certain
continuity and regularity properties on $H$ and $f$ which simplify
our analysis and ensure the existence of $\mathcal{H}f$ for every
appropriate $f$. The result we have found of most use is to
determine, for fixed sample path, under which conditions $\mathcal{H}$ takes the space of
bounded continuous functions on $[0,\infty)$ into the space of
bounded continuous functions on $[0,\infty)$ with a limit at
infinity.

This may be thought of as an analogue of the theory of admissibility
of (deterministic) linear continuous Volterra operators, especially
in the important case where the operator takes $BC_l(0,\infty)$ into
itself, or when $\mathcal{H}$ takes $BC$ into $BC_l$. Corduneanu has
done significant work on the general theory of admissibility for
Volterra integral operators (see~\cite{cord:1965} and
\cite{cord:1973}). One motivation for the development of such an
admissibility theory in the deterministic case is to give precise
asymptotic information regarding the solutions of integral, 
differential and integro-differential equations. Corduneanu~\cite{cord:1991} contains a comprehensive
survey of progress up to 1991, while further developments in this
theory are due to Cushing, Miller and others. More recently,
admissibility of continuous linear Volterra operators has been used to determine
asymptotic behaviour of a nonlinear integrodifferential equation
with infinite memory in Appleby, Gy\H{o}ri and
Reynolds~\cite{jaigdr:2011} in this journal. Parallel results are also available in
discrete time: indeed, recent results on the theory of admissibility
of Volterra operators is discrete time, together with applications
to Volterra summation equations, include Gy\H{o}ri and
Reynolds~\cite{igdr:2010} and Song and Baker~\cite{songbaker:2006}.

Reynolds \cite{dwr:2012} has established results which
characterise certain admissible pairs of spaces, as well as connecting the recent dynamical
systems literature with parallel, earlier work in the theory of linear operators.

Analogous to the deterministic setting, integrals of the form \eqref{eq:HfB} arise 
in the analysis of asymptotic growth rates of solutions 
of affine stochastic Volterra equations.
To illustrate we define some notation. For a function $x:[-\tau,\infty)\to\mathbb{R}$, we define the segment of $x$ at time $t\geq0$ as the function
\[
	x_t:[-\tau,\infty)\to\mathbb{R}: x_t(u)=x(t+u).
\]
Consider the affine stochastic functional equation
\begin{equation}\label{eq:XX}
	dX(t) = L(t,X_t) dt + \sigma(t) dB(t), \quad t\geq0; \quad X(0)=x_0\in\R,
\end{equation}
where $L$ is a linear functional and $\sigma$ is a continuous function.
The associated resolvent equation arises from setting $\sigma\equiv0$, giving
\begin{equation}\label{eq:xx}
	r'(t) = L(t,r_t), \quad t\geq0; \quad r(0)=1; \quad r(t)=0, \quad t<0.
\end{equation}
Providing both \eqref{eq:XX} and \eqref{eq:xx} have well-defined solutions then $X$ may be expressed in terms of $r$, i.e.
\begin{equation}\label{eq:Xx}
	X(t) = r(t)x_0 + \int_{0}^{t}r(t-s)\sigma(s) dB(s), \quad t\geq0.
\end{equation}
If one were to scale the solution of \eqref{eq:Xx} by a growing or decaying term then stochastic Volterra integrals of the form \eqref{eq:HfB} arise. Establishing results about the convergence of these integrals (i.e. the integrals in \eqref{eq:HfB}) therefore amounts to determining very precisely their almost sure asymptotic rates of growth or decay.
The authors propose to follow this programme of research in later works.
Another work where the scaled solution of a stochastic integral equation tends to non-trivial limit includes Appleby~\cite{app:2004gyori}.

The main results of this article concerning mean square convergence and almost sure convergence of \eqref{eq:HfB} are given in Section~\ref{sect:stochlim}. The proofs of all results are given in subsequent sections.


\subsection{Preliminaries}
Let $\R$ be the set of real numbers.
We denote by $\Rp$ the half-line $[0,\infty)$. 
 If $d$ is a positive integer, $\R^d$ is the space of
$d$-dimensional column vectors with real components and
$\R^{d_1\times d_2}$ is the space of all $d_1 \times d_2$ real
matrices. 


Let $BC(\Rp;\R^{d_1\times d_2})$ denote
the space matrices whose elements are  bounded continuous functions. The abbreviation \text{\it a.e.} stands for \text{\it almost everywhere}. The space of continuous and continuously differentiable functions on $\Rp$ with values in $\R^{d_1\times d_2}$ is denoted by $C(\Rp;\R^{d_1\times d_2})$ and $C^1(\Rp;\R^{d_1\times d_2})$ respectively. While $C^{1,0}(\Delta;\R^{d_1\times d_2})$ represents the space of functions which are continuously differentiable in their first argument and continuous in their second argument, over some two--dimensional space $\Delta$. The space of $p^{th}-$integrable functions is denoted by
\[
    L^p(\Rp;\R) := \{ f:\Rp\to\R:\int_{0}^{\infty}|f(s)|^p \,ds <+\infty \}.
\]

For any vector $x\in\R^d$ the norm $\norm{\cdot}$ denotes the Euclidean norm,
     $\norm{x}^2=\sum_{j=1}^{d}|x_j|^2$.
While for matrices, $A=(a_{i,k})\in\R^{n\times d}$, we use the Frobenius norm, i.e.
\[
    \norm{A}_{F}^2 = \sum_{i=1}^{n}\sum_{k=1}^{d}|a_{i,k}|^2.
\]
As both $\R^d$ and $\R^{n\times d}$ are finite dimensional Banach
spaces all norms are equivalent in the sense that for any other
norm,  $\norm{\cdot}$, one can find universal constants
$d_1(n,d)\leq d_2(n,d)$ such that
\[
    d_1\norm{A}_F \leq \norm{A} \leq d_2\norm{A}_F.
\]
Thus there is no loss of generality in using the Euclidean and Frobenius norms, 
which for ease of calculation, are used throughout the proofs of this paper.
Moreover we remark that the Frobenius norm is a {\it consistent
matrix norm}, i.e. for any $A\in\R^{n_1\times n_2},
B\in\R^{n_2\times n_3}$
\[
    \norm{AB}_F \leq \norm{A}_F \norm{B}_F.
\]
%

We define the following modes of convergence:
\begin{definition}
The $\mathbb{R}^{n}$-valued stochastic process $\{X(t)\}_{t\geq0}$
converges in mean-square to $X_{\infty}$ if
\[
    \lim_{t\to\infty}\mathbb{E}[\norm{X(t)-X_{\infty}}^2] = 0.
\]
\end{definition}
\begin{definition}
If there exists a $\mathbb{P}$--null set $\Omega_0$ such that for
every $\omega\not\in\Omega_0$ the following holds
\[
    \lim_{t\to\infty}X(t,\omega) = X_{\infty}(\omega),
\]
then we say $X$ converges almost surely (a.s.) to $X_{\infty}$.
\end{definition}


\section{Stochastic Limit Relation}\label{sect:stochlim}
\subsection{Mean Square Convergence}
Let $B(t)=\{B_1(t),B_2(t),...,B_d(t)\}$ be a vector of mutually
independent standard Brownian motions. 
We consider the following hypotheses: let $\Delta\subset
\mathbb{R}^2$ be defined by
$\Delta=\{(t,s): 0\leq s\leq t<+\infty\}$ 
and
\begin{gather}
\label{eq.Hcns2}
\text{$H:\Delta \to \mathbb{R}^{n\times n}$ is continuous}.
\end{gather}
We first characterise, for $f\in C([0,\infty);\mathbb{R}^{n\times
d})$ with bounded norm, the convergence of the stochastic process
$X_f=\{X_f(t):t\geq 0\}$ defined by
\begin{equation}\label{eq:Xf}
X_f(t)=\int_0^t H(t,s)f(s)\,dB(s), \quad t\geq 0
\end{equation}
to a limit as $t\to\infty$ in \textit{mean--square}.

Before discussing this convergence, we note that \eqref{eq.Hcns2} is sufficient to guarantee that 
$X_f(t)$ is a well--defined $\mathcal{F}^{B}(t)$-adapted random variable for each fixed $t$. 
Therefore the family of random variables $\{X_f(t):t\geq 0\}$ is well--defined, and $X_f$ is indeed a process. 
Condition \eqref{eq.Hcns2} also guarantees that $\mathbb{E}[X_f(t)^2]<+\infty$ for each $t\geq 0$.
Since $f\mapsto X_f$ is linear, and the family $(X_f(t))_{t\geq 0}$ is Gaussian for
each fixed $f$, the limit should also be Gaussian and linear in $f$,
as well as being an $\mathcal{F}^B(\infty)$--measurable random
variable. Therefore, a reasonably general form of the limit should
be
\begin{equation}\label{eq:Xfst}
X_f^\ast:=\int_0^\infty H_\infty(s) f(s)\,dB(s),
\end{equation}
where we would expect $H_\infty$ to be a function independent of
$f$. 
In our first main result, we show that $X_f(t)\to X_f^\ast$ in
mean square as $t\to\infty$ for each $f$.

\begin{theorem} \label{thm.msqcharacterise2}
Suppose that $H$ obeys \eqref{eq.Hcns2}.
Then the statements
\begin{itemize}
\item[(A)] There exists
$H_\infty\in C([0,\infty);\mathbb{R}^{n\times n})$ such that
$\int_{0}^{\infty}\norm{H_{\infty}(s)}^2 ds <+\infty$  and
\begin{equation}\label{eq.HtoHinfty2}
    \lim_{t\to\infty} \int_{0}^t \norm{H(t,s)-H_\infty(s)}^2\,ds=0.
\end{equation}
\item[(B)] There exists $H_\infty \in C([0,\infty);\mathbb{R}^{n\times n})$ such that
for each $f\in BC(\Rp;\mathbb{R}^{n\times d})$,
\begin{equation}\label{eq.meansquareconv2}
    \lim_{t\to\infty} \mathbb{E}\left[\norm{\int_0^t H(t,s)f(s)\,dB(s) - \int_0^\infty H_\infty(s)f(s)\,dB(s) }^2\right]=0
\end{equation}
\end{itemize}
are equivalent.
\end{theorem}

In the deterministic admissibility theory, the assumptions for
convergence are given in a different form from
\eqref{eq.HtoHinfty2}, c.f. e.g. Theorem~A.1 from \cite{jaigdr:2011}. Our next result shows that the natural
analogues of those assumptions are equivalent to
\eqref{eq.HtoHinfty2}.
\begin{proposition} \label{prop.connectdetstochmsq2}
Suppose that $H$ obeys \eqref{eq.Hcns2}. Then the following are
equivalent:
\begin{itemize}
\item[(A)] $H$ obeys \eqref{eq.HtoHinfty2};
\item[(B)] There exists $H_\infty\in C([0,\infty);\mathbb{R}^{n\times n})$ such that
\begin{gather} \label{eq.Wis02}
\lim_{T\to\infty} \limsup_{t\to\infty}\int_T^t \norm{H(t,s)}^2\,ds = 0,\\
\label{eq.HtoHinftycompact2} \lim_{t\to\infty} \int_0^T
\norm{H(t,s)-H_\infty(s)}^2 \,ds=0, \quad \text{for every $T>0$}.
\end{gather}
\end{itemize}
\end{proposition}

\subsection{Necessary Condition for Almost Sure Convergence}
We now consider the almost sure convergence of $X_f(t)$ as
$t\to\infty$ to a limit. Our next main result shows that if we have
convergence in an a.s. sense, we must also have convergence in a
mean square sense.
\begin{theorem} \label{thm.ascharacterise2}
Suppose that $H$ obeys \eqref{eq.Hcns2} and there exists $H_\infty
\in C([0,\infty);\mathbb{R}^{n\times n})$ such that for each
$f\in BC([0,\infty);\mathbb{R}^{n\times d})$,
\begin{equation}\label{eq.asconv2}
\lim_{t\to\infty} \int_0^t H(t,s)f(s)\,dB(s) = \int_0^\infty
H_\infty(s)f(s)\,dB(s), \quad \text{a.s.}
\end{equation}
Then \eqref{eq.HtoHinfty2} and \eqref{eq.meansquareconv2} hold.
\end{theorem}

Theorem~\ref{thm.msqcharacterise2} is concerned with moment behaviour of $X_f(t)=\int_0^t H(t,s)f(s)\,dB(s)$, indeed the continuity of these \emph{moments} is guaranteed by the assumption~\eqref{eq.Hcns2}.
In Theorem~\ref{thm.ascharacterise2} the condition~\eqref{eq.asconv2} may implicitly impose continuity of the \emph{sample paths} of $X_f$. The issue of continuous sample paths of  $X_f$ is addressed in Lemma~2.D. of Berger and Mizel~\cite{BergMiz:2}.
Specifically, let $H$ obey \eqref{eq.Hcns2}. Suppose that $H$ obeys a H\"older continuity condition of the following form: there exists a function $K(s)$ and a constant $\alpha>0$ such that
\[
	\int_{0}^{T}|K(s)|^2ds<+\infty
\]
and
\begin{equation}\label{eq:Holdcont}
	|H(t_2,s)-H(t_1,s)|\leq K(s)\, (t_2-t_1)^{\alpha}, \quad \text{for $0\leq s\leq t_1\leq t_2 \leq T$}.
\end{equation}
Since $H$ is continuous, it follows that there exist constants $\epsilon>0$, $D>0$ such that
\[
	\sup_{t\in[0,T]}\int_{0}^{T}|H(t,s)|^{2+\epsilon}ds \leq D.
\]
Lemma~2.D. of \cite{BergMiz:2} now guarantees that a continuous version of
\[
	\int_{0}^{t}H(t,s)f(s)dB(s)
\]
exists on $[0,T]$.
\begin{remark}
Therefore we have shown that \eqref{eq.HtoHinfty2} is a necessary
condition for a.s. convergence. It is of course natural to then ask
whether \eqref{eq.HtoHinfty2} is sufficient. We show by a simple
example that in general additional conditions are needed in order
for \eqref{eq.asconv2} to hold.
It is further noted that the assumed continuity and structure of $H$ in Examples~\ref{example.1} and~\ref{example.2} immediately gives the continuity of the sample paths of $\int_{0}^{t}H(t,s)f(s)dB(s)$, and that the sufficient condition \eqref{eq:Holdcont} is not needed.
\end{remark}
\begin{example} \label{example.1}
Suppose that $H^\sharp:[0,\infty)\to \mathbb{R}$ and
$H_\infty:[0,\infty)\to \mathbb{R}$ are continuous functions, and
define
\[
H(t,s)=H_\infty(s)+H^\sharp(t), \quad (t,s)\in \Delta.
\]
Then $H$ is continuous. Suppose also that $H_\infty\in
L^2([0,\infty);\mathbb{R})$. 
By Theorem~\ref{thm.msqcharacterise2}, it follows that
\begin{equation} \label{eq.msqexample1}
\lim_{t\to\infty} \sqrt{t}H^\sharp(t)=0
\end{equation}
is necessary and sufficient for \eqref{eq.meansquareconv2}. It is also a necessary condition for
\eqref{eq.asconv2}.
If one further supposes that $H^\sharp$ obeys
\begin{equation} \label{eq.asexample1}
\lim_{t\to\infty} \sqrt{t\log\log t} H^\sharp(t)=0,
\end{equation}
then \eqref{eq.asconv2} holds.

%

Obviously the conditions \eqref{eq.asexample1} and
\eqref{eq.msqexample1} do not coincide; in fact,
\eqref{eq.asexample1} implies \eqref{eq.msqexample1}. This provides
an example of the veracity of Theorem~\ref{thm.ascharacterise2}
which can be verified independently of the general proof of that
result.

We note also that it is very difficult to relax
\eqref{eq.asexample1} and still have the integral $\int_0^t
H(t,s)f(s)\,dB(s)$ tending to a limit a.s. as $t\to\infty$. Indeed,
there exist functions $H^\sharp$ which do not obey
\eqref{eq.asexample1}, and so must satisfy
\[
\limsup_{t\to\infty} \sqrt{t\log\log t}\, |H^\sharp(t)|>0,
\]
for which
\begin{equation} \label{eq.nolimitexample1}
\mathbb{P}\left[\lim_{t\to\infty} \int_0^t H(t,s)\,dB(s)
\quad\text{exists}\right]=0,
\end{equation}
while at the same time we still have \eqref{eq.meansquareconv2}.

A choice of $H^\sharp$ which satisfies these conditions can readily
be made. Consider a continuous function $H^\sharp$ which obeys
$H^\sharp(n)=1/\sqrt{n\log\log(n+2)}$ for all integers $n\geq 1$ but
for which $\sqrt{t}H^\sharp(t)\to 0$ as $t\to\infty$ and
$\limsup_{t\to\infty} \sqrt{t\log\log t} |H^\sharp(t)|<+\infty$.
%
\end{example}
The next example shows that sometimes the conditions which give mean
square convergence and a.s. convergence are the same.
\begin{example} \label{example.2}
Suppose that $H^\sharp:[0,\infty)\to \mathbb{R}$ and
$H_\infty:[0,\infty)\to \mathbb{R}$ are continuous functions, and
define
\[
H(t,s)=H_\infty(s)H^\sharp(t), \quad (t,s)\in \Delta.
\]
Then $H$ is continuous. Suppose also that $H_\infty\in
L^2([0,\infty);\mathbb{R})$. 
Therefore, by Theorem~\ref{thm.msqcharacterise2}, we have
\eqref{eq.meansquareconv2} if and only if
\begin{equation*}
\lim_{t\to\infty} H^\sharp(t)=1.
\end{equation*}
We know by Theorem~\ref{thm.ascharacterise2} that this condition is
also necessary for a.s. convergence.
Indeed it is also sufficient for \eqref{eq.asconv2}.
%
\end{example}

\subsection{Sufficient Conditions for Almost Sure Convergence}
We now investigate sufficient conditions for a.s. convergence for
functions $H$ which need not necessarily be of the form
\[
H(t,s)=\sum_{j=1}^n H_j(s)H_j^\sharp(t), \quad (t,s)\in \Delta,
\]
and which are covered by explicit and direct calculations in
Examples~\ref{example.1} and \ref{example.2}. 

\begin{theorem} \label{thm.assuff2}
Suppose that $H$ obeys \eqref{eq.Hcns2} and also that $H\in
C^{1,0}(\Delta;\mathbb{R}^{n\times n})$. Suppose also that there
exists $H_\infty\in C([0,\infty);\mathbb{R}^{n\times n})$ such that
$\int_{0}^{\infty}\norm{H_{\infty}(s)}^2 ds <+\infty$ and
\begin{equation} \label{eq.Htilderateto02}
\lim_{t\to\infty} \int_0^t \norm{H(t,s)-H_\infty(s)}^2\,ds \cdot
\log t=0,
\end{equation}
and
\begin{multline}
\label{eq.H1to02}
\text{There exists $q\geq 0$ and $c_q>0$ such that } \\
\int_0^t \norm{H_1(t,s)}^2\,ds \leq c_q(1+t)^{2q}, \quad
\norm{H(t,t)}^2\leq c_q(1+t)^{2q}.
\end{multline}
Then $H$ obeys \eqref{eq.asconv2}.
\end{theorem}
\begin{remark} \label{rem.strengthenedtildeHto0}
Condition \eqref{eq.Htilderateto02} implies a rate of
decay of $\int_0^t \norm{H(t,s)-H_\infty(s)}^2\,ds$ to zero as
$t\to\infty$. This strengthens the hypothesis \eqref{eq.HtoHinfty2}
which is known, by Theorem~\ref{thm.ascharacterise2}, to be
necessary.
\end{remark}

\begin{remark}
We comment now on the continuity of the sample paths of the process $\int_{0}^{t}H(t,s)f(s)dB(s)$ in Theorem~\ref{thm.assuff2}. For simplicity we discuss only the scalar case.
This continuity of the sample paths is assured by the derivative condition \eqref{eq.H1to02}. Fix $T>0$ and let $0\leq s\leq t_1\leq t_2\leq T$. Then, as $H\in C^{1,0}(\Delta;\mathbb{R})$, by the Mean Value Theorem we get 
\[
	|H(t_2,s)-H(t_1,s)| = |H_1(t^*,s)|\, |t_2-t_1|, 
\]
for some $t^*=t^\ast(s)\in[t_1,t_2]$. Next, define $K(s):=\sup_{t_1\leq t\leq t_2}|H_1(t,s)|$. This is well--defined and finite by the 
continuity of $H_1$. Therefore 
\[
|H(t_2,s)-H(t_1,s)| 
\leq K(s)\, |t_2-t_1|, \text{ for all $0\leq s\leq t_1\leq t_2\leq T$},
\]
which is \eqref{eq:Holdcont} with $\alpha=1$. Note moreover that the continuity of $s\mapsto K(s)$ ensures that 
\[
\int_0^T |K(s)|^2\,ds <+\infty,
\]
and therefore all the conditions of Lemma~2.D. of Berger and Mizel~\cite{BergMiz:2} are satisfied.
\end{remark}

\begin{remark}\label{rk:Hint}
    While the pointwise bound on $\norm{H(t,t)}$ in \eqref{eq.H1to02} may appear quite mild,
    one may prefer an integral condition to this pointwise bound as this would allow $\norm{H(t,t)}$ to potentially have ``thin spikes" of larger than polynomial order. Scrutiny of the proof of Theorem~\ref{thm.assuff2} reveals that the condition $\norm{H(t,t)}^2 \leq c_q(1+t)^{2q}$ can be replaced by
\begin{equation}\label{eq:Hint2}
	\lim_{k\to\infty}\int_{k^{\theta}}^{(k+1)^{\theta}}\norm{H(s,s)}^2 \,ds \cdot \log k = 0, \quad \text{ for }0<\theta<1/(1+2q),
\end{equation}
where the limit is taken through the integers. 
Nevertheless for simplicity we retain the condition on $\norm{H(t,t)}$ in the statement of Theorem~\ref{thm.assuff2}.
\end{remark}

\begin{example}\label{ex:logcon}
In light of Examples~\ref{example.1} and~\ref{example.2} 
(both of whose hypotheses lie outside the scope of Theorem~\ref{thm.assuff2},
in spite of the conclusion of that theorem still holding) one might query the sharpness of \eqref{eq.Htilderateto02}.
In fact, a condition of the form \eqref{eq.Htilderateto02} is, to some extent, necessary. 

Suppose, for example that $H(t,s)=\e^{-(t-s)}\sigma(s)$, $0\leq s\leq t$, where $t\mapsto\sigma(t)$ is a continuous function. 
Define the process
\[
	Y(t) = \int_{0}^{t}H(t,s)dB(s), \quad t\geq0.
\]
Then, in the context of Theorem~\ref{thm.msqcharacterise2} or~\ref{thm.ascharacterise2}, $H_\infty(s)=0$.
Moreover, if $Y$ converges to a limit almost surely we note from \eqref{eq.asconv2} that the limit should be zero, and by Theorem~\ref{thm.ascharacterise2}, condition \eqref{eq.HtoHinfty2} must hold. 
Hence, we must have that $\lim_{t\to\infty}\int_{0}^{t}H(t,s)^2\,ds=0$. It can be shown that this happens if and only if $\int_{t}^{t+1}\sigma(s)^2\,ds\to0$ as $t\to\infty$.

Suppose, moreover, that $\lim_{t\to\infty}\int_{t}^{t+1}\sigma(s)^2\,ds\cdot \log t = L\in(0,\infty]$.
Then \eqref{eq.Htilderateto02} is violated. As $Y$ is the unique adapted process satisfying the stochastic differential equation 
$dY(t)=-Y(t)+\sigma(t)dB(t)$, then it is a consequence of Theorem~4.1 of Appleby, Cheng and Rodkina~\cite{AppCheRod} that 
if $L\in(0,\infty)$ then $Y$ is bounded and does not tend to a limit or otherwise $L=+\infty$ and $Y$ is unbounded,
i.e. the conclusion of Theorem~\ref{thm.assuff2} does not hold.

Moreover, \eqref{eq.Htilderateto02} may in some sense be regarded as sharp  for $Y$ tending to a limit. From Theorem~4.1 of \cite{AppCheRod} it is known that if $Y$ tends to a limit then this limit must be zero. In Theorem~4.4 of \cite{AppCheRod} and the discussion preceding it, it is shown when $\int_{n}^{n+1}\sigma(s)^2\,ds$ is a decreasing sequence, that $Y(t)\to0$ as $t\to\infty$ a.s. implies $\int_{n}^{n+1}\sigma(s)^2\,ds\cdot\log n\to0$. This implies \eqref{eq.Htilderateto02}.
On the other hand, if $\int_{n}^{n+1}\sigma(s)^2\,ds\cdot\log n\to0$ as $n\to\infty$, we have 
by part (A) of Theorem~4.1 of \cite{AppCheRod} it follows that $Y(t)\to0$ as $t\to\infty$ a.s. Therefore, in this case, we see that \eqref{eq.Htilderateto02} is necessary and sufficient for convergence.

If $Y$ is bounded we have that $\int_{t}^{t+1}\sigma(s)^2\,ds\to0$ as $t\to\infty$.
Therefore imposing a polynomial growth bound on $\sigma$ is not restrictive.
Under this restriction and the condition $\lim_{t\to\infty}\sigma(t)^2\log t=L\in(0,\infty)$ we have that \eqref{eq.Htilderateto02}
is violated but \eqref{eq.H1to02} holds.
Hence \eqref{eq.Htilderateto02} is chiefly responsible for convergence in this case.
%
\end{example}


\section{Proofs of Admissibility Results}\label{sect:adproof}
The following proofs are given for scalar valued functions. The multi-dimensional results may be obtained by considering matrix functions component-wise (in such calculations it is often convenient to use the Frobenius norm due to It\^o's isometry).
\begin{remark} \label{rem.msqimpHinftyL2}
If \eqref{eq.meansquareconv2} holds, it is implicit that the stochastic integral
\[
\int_0^\infty H_\infty(s)f(s)\,dB(s)
\]
exists for every $f\in BC((0,\infty);\mathbb{R})$, and in particular this holds in the case when $f(s)=1$ for all $s\geq 0$. Therefore we have that
$\int_0^\infty H_\infty(s)\,dB(s)$ exists. Therefore, by the martingale convergence theorem, we have that $H_\infty\in L^2((0,\infty);\mathbb{R})$.
\end{remark}
\begin{proof}[\bf Proof of Theorem~\ref{thm.msqcharacterise2}]
It is not difficult to see using It\^o's isometry that
\begin{align}
    &\mathbb{E}\left[ \left(\int_{0}^{t}H(t,s)f(s)dB(s) - \int_{0}^{\infty}H_{\infty}(s)f(s)dB(s)\right)^{2} \right] \notag \\
    &\qquad= \int_{0}^{t}\biggl(H(t,s)-H_{\infty}(s) \biggr)^2 f(s)^2 ds
        + \int_{t}^{\infty} H_{\infty}(s)^2 f(s)^{2} ds, \label{eq:Fmsq}
\end{align}
where the independence of the elements of the Brownian vector and of
stochastic integrals over non-overlapping intervals has been used.

Firstly we show that $(\text{A})$ implies $(\text{B})$. Now as $f\in BC([0,\infty);\mathbb{R})$,
\begin{align*}
	&\mathbb{E}\left[ \left(\int_{0}^{t}H(t,s)f(s)dB(s) - \int_{0}^{\infty}H_{\infty}(s)f(s)dB(s)\right)^{2} \right] \\
	&\qquad\leq\left(\int_{0}^{t}(H(t,s)-H_{\infty}(s))^{2}ds
		+ \int_{t}^{\infty} H_{\infty}(s)^2 ds \right)\sup_{s\geq0} |f(s)|^{2}.
\end{align*}
Taking the limit as $t\to\infty$, then by hypothesis both terms on the right--hand side of the above
inequality tend to zero, and so we obtain \eqref{eq.meansquareconv2}.

Suppose to the contrary that $(\text{B})$ holds. By Remark~\ref{rem.msqimpHinftyL2}, we have that $H_\infty\in L^2([0,\infty);\mathbb{R})$.
Rearranging \eqref{eq:Fmsq} with $f(s)=1$ for all $s\geq0$,
\begin{align*}
    \int_{0}^{t}\biggl(H(t,s)-&H_{\infty}(s) \biggr)^2 ds \\
        &=\mathbb{E}\left[ \left(\int_{0}^{t}H(t,s)dB(s) - \int_{0}^{\infty}H_{\infty}(s)dB(s)\right)^{2} \right]
            - \int_{t}^{\infty}H_{\infty}(s)^2 ds,
\end{align*}
Therefore by the hypothesis of $(\text{B})$, we arrive at
\[
\limsup_{t\to\infty} \int_0^t (H(t,s)-H_\infty(s))^2\,ds =0,
\]
as required.
\end{proof}

\begin{proof}[\bf Proof of Proposition~\ref{prop.connectdetstochmsq2}]
We prove that (A) implies (B) first. To prove \eqref{eq.Wis02}, note for any $t\geq T$ we have the estimate
\begin{align*}
\int_T^t H^2(t,s)\,ds
&= \int_T^t (H(t,s)-H_\infty(s)+H_\infty(s))^2\,ds\\
&\leq 2\int_0^t (H(t,s)-H_\infty(s))^2 \,ds+  2\int_T^t H_\infty^2(s)\,ds.
\end{align*}
Since $H_\infty\in L^2([0,\infty);\mathbb{R})$ and \eqref{eq.HtoHinfty2} holds, we have
\[
\limsup_{t\to\infty} \int_T^t H^2(t,s)\,ds \leq \int_T^\infty H_\infty^2(s)\,ds.
\]
Since the lefthand side is monotone in $T$, we may take the limit as $T\to\infty$ on both sides, using
the fact that $H_\infty\in L^2([0,\infty);\mathbb{R})$ to obtain the desired conclusion.

To show \eqref{eq.HtoHinftycompact2}, let $T>0$ be arbitrary. Then, for any $t\geq T$ we have
\[
\int_0^T (H(t,s)-H_\infty(s))^2 \,ds\leq \int_0^t (H(t,s)-H_\infty(s))^2 \,ds,
\]
whence the result letting $t\to\infty$ and applying \eqref{eq.HtoHinfty2}. Thus (A) implies (B).

To prove that (B) implies (A), we first must show that $H_\infty\in L^2([0,\infty);\mathbb{R})$. We start by observing
that \eqref{eq.Wis02} is nothing other than $\lim_{T\to\infty} L(T)=0$ where
\[
L(T):=\limsup_{t\to\infty}\int_T^t H(t,s)^2\,ds.
\]
Since $L$ is non--increasing, for every $\epsilon>0$ there exists $T_0(\epsilon)>0$ such that $L(T)<\epsilon$ for all $T\geq T_0(\epsilon)$.
Now, let $T\geq T_0$. Suppose also that $t\geq T$. Then
\begin{align*}
\int_{T_0}^T H_\infty^2(s)\,ds
&\leq 2\int_{T_0}^T (H_\infty(s)-H(t,s))^2\,ds +2\int_{T_0}^T H(t,s)^2\,ds\\
&\leq 2\int_{0}^T (H_\infty(s)-H(t,s))^2\,ds +2\int_{T_0}^t H(t,s)^2\,ds.
\end{align*}
So by \eqref{eq.HtoHinftycompact2} we have
\[
\int_{T_0}^T H_\infty^2(s)\,ds \leq  2L(T_0),
\]
and since the righthand side is independent of $T$, it follows that $H_\infty\in L^2([0,\infty);\mathbb{R})$, which is one part of
\eqref{eq.HtoHinfty2}.

To prove the other part, let $t\geq T>0$. Then we have the estimate
\begin{align*}
\lefteqn{\int_0^t (H(t,s)-H_\infty(s))^2\,ds}\\
&=\int_0^T (H(t,s)-H_\infty(s))^2\,ds + \int_T^t (H(t,s)-H_\infty(s))^2\,ds\\
&\leq \int_0^T (H(t,s)-H_\infty(s))^2\,ds + 2\int_T^t H(t,s)^2\,ds + 2\int_T^t H_\infty(s)^2\,ds.
\end{align*}
Since $H_\infty\in L^2([0,\infty);\mathbb{R})$ and $H$ obeys \eqref{eq.HtoHinftycompact2}, we have
\[
\limsup_{t\to\infty} \int_0^t (H(t,s)-H_\infty(s))^2\,ds\leq
2\limsup_{t\to\infty} \int_T^t H(t,s)^2\,ds + 2\int_T^\infty H_\infty(s)^2\,ds.
\]
Now letting $T\to\infty$ on both sides of the inequality and using \eqref{eq.Wis02} proves the other part of \eqref{eq.HtoHinfty2}
\end{proof}

\begin{proof}[\bf Proof of Theorem~\ref{thm.ascharacterise2}]
Suppose for a moment that $f(s)=1$ for all $s\geq 0$. Then by \eqref{eq.asconv2} it follows that $\int_0^\infty H_\infty(s)\,dB(s)$ exists.
Therefore, by the martingale convergence theorem, we have that $H_\infty\in L^2([0,\infty);\mathbb{R})$, which is one part of \eqref{eq.meansquareconv2}. Therefore, for $f\in BC([0,\infty);\mathbb{R})$ the processes $X_f$ and $X_f^\infty$ in \eqref{eq:Xf} and \eqref{eq:Xfst} are well--defined. Also,
as $H_\infty\in L^2([0,\infty);\mathbb{R})$, we have that $X_f^\ast$ is well--defined, and thus $X_f^\infty(t)\to X_f^\ast$ as $t\to\infty$ a.s., where
\[
	X_f^\ast(t) = \int_{0}^{\infty}H_{\infty}(s) dB(s).
\]
Next, define $Y_f(t):=X_f(t)-X_f^\infty(t)$ for $t\geq 0$. Evidently, we have that $\mathbb{E}[Y_f(t)]=0$ for all $t\geq 0$. Also, we have from \eqref{eq.asconv2}
that $X_f(t)\to X_f^\ast$ as $t\to\infty$ a.s. Therefore
\[
\lim_{t\to\infty} Y_f(t)=\lim_{t\to\infty}\{X_f(t)-X_f^\infty(t)\} = \lim_{t\to\infty}\{X_f(t)-X_f^\ast+X_f^\ast-X_f^\infty(t)\}=0,
\]
almost surely. Notice also that $(Y_f(t))_{t\geq 0}$ is a Gaussian process. Since it converges a.s., it does so to a Gaussian random variable which has zero mean and zero variance, and by the argument of pp304--305 in Shiryaev~\cite{Shir}, we have that $\text{Var}[Y_f(t)]\to 0$ as $t\to\infty$. Since $\mathbb{E}[Y_f(t)]=0$, we also have $\mathbb{E}[Y_f^2(t)]\to 0$ as $t\to\infty$. However, by It\^o's isometry,
\begin{align*}
\mathbb{E}[Y_f^2(t)]&=\mathbb{E}\left[(X_f(t)-X_f^\infty(t))^2\right]
=\int_0^t (H(t,s)-H_\infty(s))^2f^2(s)\,ds.
\end{align*}
Therefore we have
\[
\lim_{t\to\infty} \int_0^t (H(t,s)-H_\infty(s))^2f^2(s)\,ds = 0,
\]
and choosing $f(s)=1$ for all $s\geq 0$, we arrive at the rest of \eqref{eq.HtoHinfty2}. 
Clearly \eqref{eq.meansquareconv2} now holds by virtue of Theorem~\ref{thm.msqcharacterise2}.
\end{proof}


\section{Proof of Theorem~\ref{thm.assuff2}}\label{sect:proofas}
Define $\tilde{H}=H-H_\infty$.  Notice that $H_\infty\in
L^2([0,\infty);\mathbb{R})$ implies that
\[
\lim_{t\to\infty}
\int_t^\infty H_\infty(s)f(s)\,dB(s)=0, \quad\text{a.s.}
\]
so that proving
\begin{equation} \label{eq.asconvHtilde}
\lim_{t\to\infty} \int_0^t \tilde{H}(t,s)f(s)\,dB(s)=0, \quad
\text{a.s.}
\end{equation}
is equivalent to establishing \eqref{eq.asconv2}.

Since $H\in C^{1,0}$, we have $\tilde{H}_1=H_1$. Therefore, we have
\[
\tilde{X}_f(t):=\int_0^t \tilde{H}(t,s)f(s)\,dB(s) = \int_0^t
\left(\tilde{H}(s,s)f(s)+\int_s^t \tilde{H}_1(u,s)f(s)\,du\right)
\,dB(s).
\]
By a stochastic Fubini theorem, \cite[Theorem~4.6.64, pp.210--211]{prot}, we have
\[
\tilde{X}_f(t) = \int_0^t \tilde{H}(s,s)f(s)\,dB(s) +\int_0^t
\left(\int_0^u H_1(u,s)f(s)\,dB(s)\right)\,du.
\]
Now, let $(t_n)_{n\geq 0}$ be an increasing sequence with $t_0=0$
and $t_n\to\infty$ as $n\to\infty$. In fact, choose
\begin{equation} \label{def.tn}
t_n=n^\theta, \quad \text{for some $\theta\in
(0,1/(1+q)\wedge1/(1+2q))\subset (0,1)$},
\end{equation}
where $q$ is the number in \eqref{eq.H1to02}.

Therefore for $t\in [t_n,t_{n+1})$, we have
\begin{multline*}
\tilde{X}_f(t) = \tilde{X}_f(t_n)+\int_{t_n}^t H(s,s)f(s)\,dB(s)-\int_{t_n}^t H_\infty(s)f(s)\,dB(s) \\
+\int_{t_n}^t \left(\int_0^u H_1(u,s)f(s)\,dB(s)\right)\,du.
\end{multline*}
Hence
\begin{multline}\label{eq.master1}
\sup_{t_n\leq t\leq t_{n+1}}|\tilde{X}_f(t)| \leq  |\tilde{X}_f(t_n)|+\sup_{t_n\leq t\leq t_{n+1}}\left|\int_{t_n}^t H_\infty(s)f(s)\,dB(s)\right| \\
+\sup_{t_n\leq t\leq t_{n+1}}\left|\int_{t_n}^t
H(s,s)f(s)\,dB(s)\right| +\int_{t_n}^{t_{n+1}} \left|\int_0^u
H_1(u,s)f(s)\,dB(s)\right|\,du.
\end{multline}
We now show that each of the four terms on the righthand side of
\eqref{eq.master1} tends to zero as $n\to\infty$ a.s.

\textbf{STEP 1: First term on the righthand side of
\eqref{eq.master1}.} First we prove that
\begin{equation} \label{eq.Xtnto0}
\lim_{n\to\infty} \tilde{X}_f(t_n)=0, \quad\text{a.s.}
\end{equation}
Notice that $\tilde{X}_f(t_n)$ is normally distributed with mean
zero and variance $v_n^2$ where
\[
v_n^2:=\int_0^{t_n} \tilde{H}^2(t_n,s)f^2(s)\,ds\leq \int_0^{t_n}
\tilde{H}^2(t_n,s)\,ds \cdot \sup_{s\geq 0}f^2(s).
\]
Using \eqref{eq.Htilderateto02} and the fact that $t_n\to\infty$ as
$n\to\infty$, we have
\[
\lim_{n\to\infty} \int_0^{t_n} \tilde{H}(t_n,s)^2\,ds\cdot \log
t_n=0,
\]
Therefore
\begin{equation} \label{eq.vnto0rate}
\limsup_{n\to\infty}  v_n^2 \log t_n \leq \limsup_{n\to\infty}
\int_0^{t_n} \tilde{H}(t_n,s)^2\,ds \cdot \sup_{s\geq 0}f^2(s)\cdot
\log t_n =0.
\end{equation}
Since $X_n:=\tilde{X}_f(t_n)/v_n$ is a standardised normal random
variable, we have that
\[
\limsup_{n\to\infty} \frac{|\tilde{X}_f(t_n)|}{\sqrt{2} v_n(\log
n)^{1/2}} = \limsup_{n\to\infty} \frac{|X_n|}{\sqrt{2\log n}}\leq 1,
\quad \text{a.s.},
\]
the last inequality being a routine consequence of the
Borel--Cantelli lemma. Thus
\begin{align*}
\limsup_{n\to\infty} |\tilde{X}_f(t_n)| &=
\limsup_{n\to\infty} \frac{|\tilde{X}_f(t_n)|}{\sqrt{2} v_n(\log n)^{1/2}} \cdot \sqrt{2} v_n(\log n)^{1/2}\\
&\leq \sqrt{2} \limsup_{n\to\infty} v_n (\log t_n)^{1/2}
\sqrt{\frac{\log n}{\log t_n}}=0,
\end{align*}
due to \eqref{def.tn} and \eqref{eq.vnto0rate}, proving
\eqref{eq.Xtnto0}.

\textbf{STEP 2: Second term on the righthand side of
\eqref{eq.master1}} Next we show that
\begin{equation} \label{eq.HinftyInto0}
\lim_{n\to\infty} \sup_{t_n\leq t\leq t_{n+1}}\left| \int_{t_n}^t
H_\infty(s)f(s)\,dB(s)\right|=0, \quad\text{a.s.}
\end{equation}
To do this, notice for every $\epsilon>0$ by Chebyshev's inequality
and the Birkholder--Davis--Gundy inequality, c.f. e.g. \cite[Theorem~1.3.8, Theorem~1.7.3]{Mao2} that
\begin{align*}
\lefteqn{\mathbb{P}\left[\sup_{t_n\leq t\leq t_{n+1}}\left| \int_{t_n}^t H_\infty(s)f(s)\,dB(s)\right|>\epsilon\right]}\\
&\leq \frac{1}{\epsilon^2}\mathbb{E}\left[\sup_{t_n\leq t\leq t_{n+1}}\left| \int_{t_n}^t H_\infty(s)f(s)\,dB(s)\right|^2\right]\\
&\leq \frac{4}{\epsilon^2} \mathbb{E}\left[\left|\int_{t_n}^{t_{n+1}} H_\infty(s)f(s)\,dB(s)\right|^2\right]\\
&= \frac{4}{\epsilon^2} \int_{t_n}^{t_{n+1}}
H_\infty^2(s)f^2(s)\,ds.
\end{align*}
Since $H_\infty\in L^2([0,\infty);\mathbb{R})$, and $f\in
BC([0,\infty);\mathbb{R})$, we have
\[
\sum_{n=0}^\infty \mathbb{P}\left[\sup_{t_n\leq t\leq t_{n+1}}\left|
\int_{t_n}^t H_\infty(s)f(s)\,dB(s)\right|>\epsilon\right] \leq
\frac{4}{\epsilon^2} \int_0^\infty H_\infty^2(s)f^2(s)\,ds.
\]
By the Borel--Cantelli Lemma, we have that \eqref{eq.HinftyInto0}
holds.

\textbf{STEP 3: Third term on the righthand side of
\eqref{eq.master1}.}
\[
\lim_{n\to\infty} U_n=0, \quad \text{a.s.}
\]
where
\[
U_n=\sup_{t_n\leq t\leq t_{n+1}} \left| \int_{t_n}^t
H(s,s)f(s)\,dB(s)\right|.
\]
Note that $(U_n)_{n\geq 0}$ is a sequence of independent random
variables.

Notice that on the interval $[t_n,t_{n+1}]$, by the martingale time
change theorem, there exists a Brownian motion $\tilde{B}$ such that
\begin{align*}
U_n
&= \sup_{t_n\leq t\leq t_{n+1}}\left|\tilde{B}_n\left(\int_{t_n}^t H^2(s,s)f^2(s)\,ds\right)\right|\\
&=\sup_{0\leq \tau\leq \int_{t_n}^{t_{n+1}} H^2(s,s)f^2(s)\,ds}|\tilde{B}_n(\tau)|\\
&\leq \sup_{0\leq \tau\leq \int_{t_n}^{t_{n+1}} H^2(s,s)\,ds\cdot
\sup_{v\geq 0}f^2(v)}|\tilde{B}_n(\tau)|.
\end{align*}
Therefore, with $w_n:=\int_{t_n}^{t_{n+1}} H^2(s,s)\,ds\cdot
\sup_{v\geq 0}f^2(v)$, we have for some Brownian motion $W$ that
\[
\mathbb{P}[U_n>\epsilon]\leq \mathbb{P}[\sup_{0\leq \tau\leq
w_n}|W(\tau)|>\epsilon].
\]
Using the symmetry of the distribution function leads to the estimate
\[
\mathbb{P}[U_n>\epsilon]
\leq 2 \mathbb{P}[|W(w_n)|>\epsilon]
\leq 4\mathbb{P}[W(w_n)>\epsilon]=4\mathbb{P}[Z>\epsilon/\sqrt{w_n}],
\]
where $Z$ is a standard normal random variable, and we  interpret
the right hand side as zero if $w_n=0$. Hence if $\Phi$ is the
distribution function of a standard normal random variable and
\[
\sum_{n=0}^\infty
\left\{ 1-\Phi\left(\frac{\varepsilon}{\sqrt{\int_{t_n}^{t_{n+1}}
H^2(s,s)\,ds}}\right)\right\}<+\infty, \quad \text{for all $\varepsilon>0$},
\]
we have that $\lim_{n\to\infty} U_n=0$, a.s. The sum is finite
provided
\[
\lim_{n\to\infty} \int_{t_n}^{t_{n+1}} H^2(s,s)\,ds\cdot \log n=0.
\]
Since $H(t,t)^2\leq c_q(1+t^{2q})$, we have that
\[
\int_{t_n}^{t_{n+1}} H^2(s,s)\,ds\cdot \log n \leq  c_q
\int_{n^\theta}^{(n+1)^\theta} \{1+s^{2q}\}\,ds \cdot \log n,
\]
so the right hand side is of the order
$n^{-1+\theta}n^{2q\theta}\log n=n^{-1+(2q+1)\theta}\log n\to 0$ as
$n\to\infty$, because $\theta<1/(1+2q)$.

\textbf{STEP 4: Fourth term on the righthand side of
\eqref{eq.master1}.} Finally, we show that
\begin{equation} \label{eq.fourthterm}
\lim_{n\to\infty} Z_n=0,\quad \text{a.s.}
\end{equation}
where
\begin{equation} \label{def.Zn}
Z_n:=\int_{t_n}^{t_{n+1}} \left|\int_0^u
H_1(u,s)f(s)\,dB(s)\right|\,du.
\end{equation}
By \eqref{eq.H1to02} there exists $c_q>0$ such that
\[
\int_0^t H_1^2(t,s)\,ds \leq c_q(1+t)^{2q}, \quad t\geq 0.
\]
By \eqref{def.tn}, $\theta<1/(1+q)\leq 1$, so we can choose $p\in
\mathbb{N}$ so large that $2p[1-(1+q)\theta]>1$. Clearly for such a
$p\in \mathbb{N}$ we have, via Jensen's inequality
\[
Z_n^{2p}\leq (t_{n+1}-t_n)^{2p-1}\int_{t_n}^{t_{n+1}} \left(\int_0^u
H_1(u,s)f(s)\,dB(s)\right)^{2p}\,du,
\]
so there exists $C_p>0$ such that
\begin{align*}
\mathbb{E}[Z_n^{2p}]&\leq C_p(t_{n+1}-t_n)^{2p-1}\int_{t_n}^{t_{n+1}} \left(\int_0^u H_1^2(u,s)f^2(s)\,ds\right)^{p}\,du\\
&\leq C_p\sup_{s\geq
0}f^{2p}(s)(t_{n+1}-t_n)^{2p-1}\int_{t_n}^{t_{n+1}} \left(\int_0^u
H_1^2(u,s)\,ds\right)^{p}\,du.
\end{align*}
Then
\begin{align*}
\mathbb{E}[Z_n^{2p}]
&\leq C_p\sup_{s\geq 0}f^{2p}(s)(t_{n+1}-t_n)^{2p-1}\int_{t_n}^{t_{n+1}} \left(c_q(1+u)^{2q}\right)^{p}\,du\\
&\leq C_p c_q^p\sup_{s\geq 0}f^{2p}(s) \cdot (t_{n+1}-t_n)^{2p-1}\int_{t_n}^{t_{n+1}} (1+u)^{2qp}\,du\\
&\leq C_p c_q^p\sup_{s\geq 0}f^{2p}(s) \cdot (t_{n+1}-t_n)^{2p}
(1+t_{n+1})^{2qp}.
\end{align*}
Since $t_n=n^\theta$, the right hand side is of the order
$[n^{\theta-1}]^{2p} n^{2pq\theta}=n^{-2p[1-(1+q)\theta]}$ as
$n\to\infty$. By Chebyshev's inequality, for any $\epsilon>0$ we
have
\[
\mathbb{P}[|Z_n|>\epsilon]\leq
\frac{1}{\epsilon^{2p}}\mathbb{E}[Z_n^{2p}] \leq C_{\epsilon,p}
n^{-2p[1-(1+q)\theta]},
\]
and because $2p[1-(1+q)\theta]>1$, the righthand side is summable.
Therefore, by the Borel--Cantelli lemma, we have
\eqref{eq.fourthterm}.


\subsection{Proof of Examples}

\begin{proof}[Proof of Example~\ref{example.1}]
Since $H_\infty$ is in $L^2([0,\infty);\mathbb{R})$, for each $f\in
BC([0,\infty);\mathbb{R})$ we have that the integral $\int_0^t
H_\infty(s)f(s)\,dB(s)$  tends almost surely as $t\to\infty$ to
$\int_0^\infty H_\infty(s)f(s)\,dB(s)$. 
Thus
\begin{multline}  \label{eq.limitexample1}
\int_0^t H(t,s)f(s)\,dB(s)-\int_0^\infty H_\infty(s)f(s)\,dB(s)
\\=H^\sharp(t)\int_0^t f(s)\,dB(s) - \int_t^\infty H_\infty(s)f(s)\,dB(s).
\end{multline}
As $H^\sharp$ obeys \eqref{eq.asexample1}
we have $H^\sharp(t)\to 0$ as $t\to\infty$. If $f\in
L^2([0,\infty);\mathbb{R})$, then $\int_0^t f(s)\,dB(s)$ tends to a
finite limit a.s., and therefore both terms on the righthand side of
\eqref{eq.limitexample1} tend to zero as $t\to\infty$ a.s., and
\eqref{eq.asconv2} holds.

On the other hand, if $f\not\in L^2([0,\infty);\mathbb{R})$, then
the martingale time change theorem and the Law of the Iterated
Logarithm, c.f. e.g \cite[Exercise ~5.1.15]{RevYor} 
as well as the boundedness of $f$ and \eqref{eq.asexample1} give that
$H^\sharp(t)\int_0^t f(s)\,dB(s)\to 0$ as $t\to\infty$ a.s. Hence we
have that \eqref{eq.asconv2} holds.

We turn now to the question of relaxing \eqref{eq.asexample1} and still having \eqref{eq.asconv2} holding.
By virtue of the fact that $H^\sharp$ obeys \eqref{eq.msqexample1},
we have that \eqref{eq.meansquareconv2} holds. By the Law of the
Iterated Logarithm, we have that
\[
\limsup_{t\to\infty}|H^\sharp(t)B(t)|<+\infty, \quad \text{a.s.}
\]
However,
\[
\limsup_{t\to\infty}|H^\sharp(t)B(t)|\geq
\limsup_{n\to\infty}|H^\sharp(n) B(n)| =\sqrt{2}\limsup_{n\to\infty}
\frac{|B(n)|}{\sqrt{2n\log\log(n+2)}}=\sqrt{2},
\]
a.s., by the discrete version of the Law of the iterated logarithm, cf. e.g. \cite[Theorem~10.2.1]{chowt}.
If $f(t)=1$ for all $t\geq 0$ in \eqref{eq.limitexample1}, we have
\[
\int_0^t H(t,s)\,dB(s)-\int_0^\infty
H_\infty(s)\,dB(s)=H^\sharp(t)B(t) - \int_t^\infty
H_\infty(s)\,dB(s).
\]
The second term on the righthand side has zero limit as $t\to\infty$
a.s., but by the above argument, the first term obeys
\[
0<\limsup_{t\to\infty}|H^\sharp(t)B(t)|<+\infty, \quad\text{a.s.}
\]
and therefore \eqref{eq.nolimitexample1} holds, as claimed.
\end{proof}

\begin{proof}[Proof of Example~\ref{example.2}]
To show that it is sufficient, suppose $f\in
BC([0,\infty);\mathbb{R})$. 
Therefore, we have the identity
\begin{multline*}
\int_0^t H(t,s)f(s)\,dB(s)-\int_0^\infty H_\infty(s)f(s)\,dB(s)\\
=(H^\sharp(t)-1)\int_0^t H_\infty(s)f(s)\,dB(s)-\int_t^\infty
H_\infty(s)f(s)\,dB(s).
\end{multline*}
Since $H^\sharp(t)\to 1$ as $t\to\infty$, the limit as $t\to\infty$
of the righthand side is zero, and so we have \eqref{eq.asconv2}.
Therefore, the condition $H^\sharp(t)\to 1$ as $t\to\infty$ is
necessary and sufficient both for \eqref{eq.asconv2} and for
\eqref{eq.meansquareconv2}.
\end{proof}

\begin{proof}[Proof of Example~\ref{ex:logcon}]
Suppose that $\lim_{t\to\infty}\int_{t}^{t+1}\sigma(s)^2\,ds=0$. Then clearly one has $\int_{n}^{n+1}\sigma(s)^2\,ds\to0$ as $n\to\infty$. Now, for $n\leq t< n+1$, where $n=\lfloor t\rfloor$,
\begin{align*}
\mathbb{E}[Y(t)^2] &= \sum_{j=1}^{n}\int_{j-1}^{j}\e^{-2(t-s)}\sigma(s)^2\,ds + \int_{n}^{t}\e^{-2(t-s)}\sigma(s)^2\,ds \\
&\leq \sum_{j=1}^{n}\e^{-2(n-j)}\int_{j-1}^{j}\sigma(s)^2\,ds + \int_{n}^{n+1}\sigma(s)^2\,ds.
\end{align*}
The first term on the right-hand side of the inequality is the convolution of a summable sequence with a sequence tending to zero and thus the convolution itself tends to zero as $t\to\infty$. Thus we have $\lim_{t\to\infty}\mathbb{E}[Y(t)^2]=0$.

Conversely, suppose that $\lim_{t\to\infty}\mathbb{E}[Y(t)^2]=0$. Then defining
\[
	y(t) = \mathbb{E}[Y(t)^2] = \int_{0}^{t}\e^{-2(t-s)}\sigma(s)^2\,ds,
\]
one obtains the identity
\[
	y(t+1) - y(t) = -2\int_{t}^{t+1}y(s)ds + \int_{t}^{t+1}\sigma(s)^2\,ds.
\]
Thus, $y(t)\to0$ as $t\to\infty$ implies $\lim_{t\to\infty}\int_{t}^{t+1}\sigma(s)^2\,ds=0$. 
Hence $\int_{t}^{t+1}\sigma(s)^2\,ds\to0$ as $t\to\infty$ completely characterises the mean square convergence of $Y$ to zero.

We now show that $\lim_{t\to\infty}\int_{t}^{t+1}\sigma(s)^2\,ds\cdot \log t = L\in(0,\infty]$ implies that
\begin{equation}\label{eq:H2log}
\liminf_{t\to\infty}\int_{0}^{t}H(t,s)^2\,ds\cdot\log t\in(0,\infty]
\end{equation}
Observe, for $n\leq t<n+1$,
\begin{align*}
\int_{0}^{t}\e^{-2(t-s)}\sigma(s)^2\,ds\cdot \log t 
&= \left( \sum_{j=1}^{n}\e^{-2t}\int_{j-1}^{j}\e^{2s}\sigma(s)^2 ds + \int_{n}^{t}\e^{-2(t-s)}\sigma(s)^2 ds \right) \log t \\
&\geq\e^{-4}\frac{\sum_{j=1}^{n}\frac{\e^{2j}}{\log(j+1)}\int_{j-1}^{j}\sigma(s)^2ds\cdot\log(j+1)}{\e^{2n}/\log n}.
\end{align*}
Standard application of Toeplitz's Lemma, \cite[Lemma~IV.3.1]{Shir}, then gives \eqref{eq:H2log}.

We show now that $\lim_{t\to\infty}\int_{t}^{t+1}\sigma(s)^2\,ds\cdot\log t=0$ implies that \eqref{eq.Htilderateto02} holds.
Observe that for $n\leq t<n+1$,
\begin{align*}
\int_{0}^{t}H(t,s)^2\,ds\cdot\log t 
&= \sum_{j=1}^{n}\int_{j-1}^{j}\e^{-2(t-s)}\sigma(s)^2 ds\log t + \int_{n}^{t}\e^{-2(t-s)}\sigma(s)^2 ds\log t \\
&\leq \frac{\sum_{j=1}^{n}\frac{\e^{2j}}{\log j}\int_{j-1}^{j}\sigma(s)^2 ds\log j}{\e^{2n}/\log(n+1)} 
+ \int_{n}^{n+1}\sigma(s)^2ds\log(n+1).
\end{align*}
Again, standard application of \cite[Lemma~IV.3.1]{Shir} gives the desired result.

Observe that $H_1(t,s) =-H(t,s)$, thus with the supposition $\int_{t}^{t+1}\sigma(s)^2ds\to0$, this gives the first part of \eqref{eq.H1to02}. The remainder of \eqref{eq.H1to02} is clearly satisfied if $\lim_{t\to\infty}\sigma(s)^2\log t =L\in(0,\infty)$.
It is obvious that $\lim_{t\to\infty}\sigma(s)^2\log t =L\in(0,\infty)$ implies $\lim_{t\to\infty}\int_{t}^{t+1}\sigma(s)^2\,ds\cdot\log t=L$ which hence violates \eqref{eq.Htilderateto02}.
\end{proof}


\end{document}